\documentclass[12pt,a4paper,twoside,final,notitlepage, leqno]{article}
\usepackage[english]{babel}
\usepackage[T1]{fontenc}
\usepackage{epsfig, graphicx, amssymb}
\usepackage{amsmath,amsthm,epsfig,amsfonts}
\usepackage{float}
\usepackage{color}
\usepackage{amssymb}

\usepackage{bm}
\usepackage{amsmath}

\usepackage{mathrsfs}
\usepackage{amsfonts}
\usepackage{latexsym}
\usepackage{multirow}
\usepackage{comment}

\newcommand{\tensorfont}[1]{\bm{\mathsf{#1}}}

\newcommand{\nueff}{\nu_{\mbox{\scriptsize eff}}}

\newcommand{\monitem}{ \smallskip \noindent $\bullet$ \quad  }
\newcommand{\moneq}{\vspace*{-7pt} \begin{equation} \displaystyle }
\newcommand{\moneqstar}{\vspace*{-6pt} \begin{equation*} \displaystyle }
\newcommand{\monendstar}{\vspace*{-6pt} \end{equation*}   }
\newcommand{\monend}{\vspace*{-7pt} \end{equation}   }
\newcommand{\moneqarraystar}{ \begin{eqnarray*} \displaystyle }
\newcommand{\monendarraystar}{ \end{eqnarray*}   }

\definecolor{vertfonce}{rgb}{0.0, 0.5, 0.0}

\def\section*#1{}

\usepackage{fancyhdr}
\renewcommand{\headrulewidth}{0pt}

\begin{document}

\fancypagestyle{plain}{ \fancyfoot{} \renewcommand{\footrulewidth}{0pt}}
\fancypagestyle{plain}{ \fancyhead{} \renewcommand{\headrulewidth}{0pt}}

~

  \vskip 2.1 cm

\centerline {\bf \LARGE Theory of the lattice Boltzmann method: }

\smallskip \smallskip

\centerline {\bf \LARGE   discrete effects due to advection }

 \bigskip  \bigskip \bigskip

\centerline { \large   Pierre Lallemand$^{a}$,  Fran\c{c}ois Dubois$^{bc}$ and Li-Shi Luo$^{ad}$}

\smallskip  \bigskip

\centerline { \it  \small
  $^a$ Beijing Computational Science Research Center, Haidian District, Beijing 100094,  China.}

\centerline { \it  \small
  $^b$   Laboratoire de Math\'ematiques d'Orsay, Facult\'e des Sciences d'Orsay,}

\centerline { \it  \small   Universit\'e Paris-Saclay, France.}

\centerline { \it  \small
$^c$    Conservatoire National des Arts et M\'etiers, LMSSC laboratory,  Paris, France.}

\centerline { \it  \small
 $^d$ Department of Mathematics \& Statistics, Old Dominion University, Norfolk, VA 23529, USA}


\bigskip  \bigskip

\centerline {27 August 2022  
{\footnote {\rm  \small $\,$ A preliminary version of this contribution was   presented 
by Pierre Lallemand
at the International Conference for Mesoscopic Methods in Engineering and Science, Hambourg (Germany), 18-22 July 2016.}}}




\bigskip  \bigskip
\noindent {\bf \large Abstract}

\noindent
Lattice Boltzmann models are briefly introduced together with
references to methods used to predict their ability for simulations of
systems described by partial differential equations that are first
order in time and low order in space derivatives.  Several previous
works have been devoted to analyzing the accuracy of these models with
special emphasis on deviations from pure Newtonian viscous behaviour,
related to higher order space derivatives of even order.  The present
contribution concentrates on possible inaccuracies of the advection
behaviour linked to space derivatives of odd order.
Detailed properties of advection-diffusion and athermal fluids are
presented for two-dimensional situations allowing to propose
situations that are accurate to third order in space
derivatives. Simulations of the advection of a gaussian dot or vortex
are presented. Similar results are discussed in appendices for
three-dimensional advection-diffusion.

\newpage 
\bigskip \bigskip    \noindent {\bf \large    1) \quad  Introduction} 

\fancyhead[EC]{\sc{Pierre Lallemand, Fran\c{c}ois Dubois and Li-Shi Luo}}
\fancyhead[OC]{\sc{Discrete effects due to advection in the lattice Boltzmann method}}
\fancyfoot[C]{\oldstylenums{\thepage}}


\label{sec:intro}

Lattice Boltzmann models have been developed over almost three
decades \cite{DDH92} based on microscopic physical models \cite{LL00}
and practices of numerical methods to solve PDE's \cite{DL09}.
The physical base is the notion of particles undergoing successive
phases of free travel and collisions. Every function of the
microscopic properties that is conserved in collisions will correspond
to a macroscopic quantity that varies slowly in space and time and
thus can be useful for computer simulations.
The kinetic theory of gases has developed relationships between elementary motions
and collisions of particles and partial differential equations describing the
behavior of the relevant macroscopic quantities. It gives guidance to setting-up
simplified models that may lead to useful numerical tools.

Computational fluid dynamics aims to predict the behavior of these
quantities.  It usually limits the description to a number of
locations in space and at a number of times. Here we choose 
$t = n \delta t$, 
$ \bm{r} := (i \hat{\bm{e}}_1 + j \hat{\bm{e}}_2) \, \delta r$ for 2-D
problems and 
$ \bm{r} := (i \hat{\bm{e}}_1 + j \hat{\bm{e}}_2 + l \hat{\bm{e}}_3)
\, \delta r$ 
for 3-D problems, which are the ``nodes'' where the state of the fluid
is defined, and $(\hat{\bm{e}}_1,\, \hat{\bm{e}}_2,\, \hat{\bm{e}}_3)$
are unit vectors for the spatial mesh. 
For simplicity further detailed  expressions will be written for the
2-D case and some results will be given for 3-D cases.

In the basic Lattice Boltzmann Model (LBM), particles move
synchronously between the various nodes, usually going to close
neighbors in one time step. This allows to define a set of $N$
elementary velocities of amplitude of the order of $\delta r/\delta
t$, $\{ c_p| p = 0,\, 1,\, \ldots,\, N-1\}$, of Cartesian components
$(c_{px},\, c_{py})$.
At time $n$, the system is fully described by a set of $N\times M$
distribution functions $f_p(i,j,n)$ or by a point $X$ is phase space
${\Phi} \in \mathbb{R}^{N\times M}$ for $M$ active nodes.
The dynamics is inspired from the Boltzmann equation. It consists in two steps:

\noindent {\it (i)} \quad $\,$ Local collision: \quad 
$f_p(i,\, j,\, n) \mapsto f_p^*(i,\, j,\, n)$

\noindent {\it (ii)} \quad Propagation to neighboring nodes: 

\quad 
$f_p(i+c_{px} , \, j+c_{py},\, n+1) = f_p^*(i,\, j,\, n)$

or from neighboring nodes:

\quad
$f_p(i, \, j,\, n+1) = f_p^*(i-c_{px},\, j-c_{py},\, n)$.

In the following, we define various Lattice Boltzmann models (Section~2),
then explain the algorithm of generationg the equivalent equations (Section~3),
and the stability analysis in the linear case (Section~4). 
Then we present 
analytic results from the linear analysis  (Section~5), 
including athermal fluid is simulated with the D2Q9  and the D2Q13 schemes.
We study the distortion of a Gaussian dot or vortex in Section~6.
Some technical precisions are presented in the appendices.


\begingroup
\let\clearpage\relax


\bigskip \bigskip    \noindent {\bf \large    2) \quad  A brief description of the lattice Boltzmann equation} 

\label{sec:LBM}

The lattice Boltzmann equation (LBE) evolves on a $d$ dimensional
lattice $\delta r \mathbb{Z}_d$ with lattice spacing $\delta r$ and is
fully defined by two ingredients: a set of discrete velocities
$\mathbb{V} := \{ \bm{c}_p \}$ and a collision model. Since the LBE is
designed to simulate low-Mach-number flows, the discrete velocities
$\mathbb{V}$ is symmetric, that is, $-\mathbb{V} = \mathbb{V}$, or,
$$
\forall\ \bm{c}_p \in \mathbb{V},
\quad
\bm{c}_{\bar{p}} := - \bm{c}_p \in \mathbb{V} ,
$$
thus, $\sum_p \bm{c}_p = \bm{0}$.
Corresponding to each discrete velocity $\bm{c}_p$, there is a
distribution function $f_p(\bm{r}_j,\, t_n)$ at every lattice point
$\bm{r}_j$ and each discrete time $t_n := n \delta t$, where $n \in
\mathbb{N}_0 := \{0,\, 1,\, 2,\, \ldots,\, \}$ and $\delta t$
is the time step size. In this setting, the unit of the velocity is $c
:= \delta r / \delta t$. The discrete velocity set $\mathbb{V}$, 
the set of nodes $\delta r \mathbb{Z}_d$, and the discrete time step size are tied
together as follows~:
$$
 \forall \, \bm{c}_p \in \mathbb{V}
 \mbox{\ and\ } \bm{r}_j \in \delta r \mathbb{Z}_d,
 \quad
 \bm{r}_j + \bm{c}_p \delta t \in \delta r \mathbb{Z}_d .
$$

The evolution of the lattice Boltzmann equation consists of two steps:
(a) a local collision model
$$
  f_p(\bm{r}_j,\, n) \mapsto f_p^*(\bm{r}_j,\, n) ,
$$
where $f_p(\bm{r}_j,\, n)$ and $f_p^*(\bm{r}_j,\, n)$ are the
pre-collision and the post-collision states at the lattice node
$\bm{r}_j$ and the time $t_n$, respectively; and (b) propagation (or
advection) from one lattice node $\bm{r}_j$ to another $\bm{r}_j +
\bm{c}_p \delta t$ in one time step according to discrete velocities
$\bm{c}_p$:
$$
f_p(\bm{r}_j + \bm{c}_p \delta t,\, n+1) = f_p^*(\bm{r}_j,\, n) .
$$
%

%

A LBM model is fully defined by two pieces of information: the set of
elementary velocities and the rules that govern the collision step.
As one usually aims to simulate fluid flows, it is highly suggested to
use a set of elementary velocities as isotropic as possible. This
means using orthogonal coordinates and for each possible velocity
amplitude, sets obtained by symmetry and permutation of the axis.

\smallskip \noindent Note that one can also use 6 velocities based on
the hexagon, but this cannot be extended to 3-D cases.

We will adopt the notation of D$d$Q$q$ for a model in $d$-dimensional
space with $q$ velocities.  In this work we shall mostly focus on the
lattice Boltzmann (LB) models in space of two dimensions (2D).
The most often used thirteen discrete velocities in 2D are listed in
Table~\ref{tab:D2Q13v}. We note that these discrete velocities conform
with the Cartesian square lattice in 2D. However, it is possible also
to use a triangular lattice in 2D \cite{DL13}.
Obviously, the Cartesian lattice in 2D can be easily extended to 3D.

\begin{table}[htbp!]    
\centering
\begin{tabular}{| r | r | l |} 
\hline
\multicolumn{1}{|c|}{Number} & 
\multicolumn{1}{c|}{$|\bm{c}_p/c|^2$} & 
\multicolumn{1}{c|}{$\bm{c}_p/c$}  
\\
\hline
 1  & 0  & $(0,\, 0)$ 
\\
 4  & 1  & $(1,\, 0)$, $(0,\, 1)$, $(-1,\, 0)$, $(0,\, -1)$ 
\\
 4  & 2  & $(1,\, 1)$, $(-1,\, 1)$, $(-1,\, -1)$, $(1,\, -1)$ 
\\
 4 & 4 & $(2,\, 0)$, $(0,\, 2)$, $(-2,\, 0)$, $(0,\, -2)$ 
\\
\hline  \end{tabular}
\caption{The first 13 discrete velocities used in the various lattice
  Boltzmann models.}
\label{tab:D2Q13v}
\end{table}

We consider simple \emph{local} collision model that gives prevalence
to the notions of ``conservation'' and symmetry, two equivalence
concepts according to N\"{o}ther. In this work we will use the linear
relaxation model proposed by d'Humi\`eres \cite{DDH92}, in which the
collision process is modeled by the linear relaxations of the velocity
moments $\{ m_p \}$ of the distribution functions $\{ f_p \}$.  Given
a set $\{ \bm{c}_p | p = 0,\, 1,\, \ldots,\, (q-1) \}$ of $q$ discrete
velocities, there always exists a $q \times q$ invertible matrix
$\tensorfont{M}$ such that
\begin{equation}
  \mathbf{m} = \tensorfont{M} \mathbf{f} ,
  \quad
  \mathbf{f} = \tensorfont{M}^{-1} \mathbf{m} ,
\end{equation}
where $\mathbf{f}$ and $\mathbf{m}$ denote the vectors of $q$
dimensions of the distribution functions $\{ f_p \}$ and the moments
$\{ m_p \}$, respectively, \emph{i.e.},
\begin{align*}
  &
  \mathbf{f}
  :=
  (f_0,\, f_1,\, \ldots,\, f_{q-1} )^\dagger ,
  \\
  &
  \mathbf{m}
  :=
  (m_0,\, m_1,\, \ldots,\, m_{q-1} )^\dagger ,
\end{align*}
where $\dagger$ denotes transpose.
It is convenient to use orthogonal polynomials on the discrete
velocity set $\mathbb{V}$ so that the relaxation processes of moments
are independent to each other.
The orthogonal polynomials with respect to a weight of unity for the models
up to thirteen velocities in 2D are given in
Table~\ref{tab:D2Q13-polynomials}.
Denote the polynomials in Table~\ref{tab:D2Q13-polynomials} by
$P_q(\bm{c}_p)$, then the transformation matrix can be constructed with
its matrix elements given by $P_q(\bm{c}_p)$, \emph{i.e.},
$\tensorfont{M}_{pq} = P_q(\bm{c}_p)$.

in the use at each node of a linear transformation of the set of
distribution functions $f_p$ to moments based on polynomials of the
elementary velocities components of increasing order chosen as
isotropic as possible. It is also convenient to orthogonalize the
moments of the same symmetry. This allows to define a ``moment
matrix'' $M$ that relates the distributions $f_p$ and the moments
$m_p$ by $m = {M} f$. (Note that $M$ must be invertible.)
We use for a 2-D model with $N$ velocities the nomenclature D2QN.
The polynomials used to generate ${M}$ (by replacing $(x,y)$ by
$(c_{px},c_{py})$ for each elementary velocity) are:


\begin{table}[htbp!]     
\centering
\begin{tabular}{|c|l|}
\hline 
model & \multicolumn{1}{c|}{Orthogonal Polynomials on $\mathbb{V}$, 
$r := \sqrt{x^2+y^2}$ }  
\\
\hline
D2Q1 & 1
\\
\hline
D2Q5 & 1, $x$, $y$, 
$-4 + 5 r^2$,
$x^2-y^2$ 
\\
\hline
D2Q9 & 1,  $x$,  $y$,  
$-4+3r^2$,
$x^2-y^2$,  $xy$,
\\
 & $-(5 - 3 r^2) x$, $-(5 - 3 r^2) y$,
$4 - \frac{3}{2} ( 7 + 3 r^2 ) r^2 $
\\
\hline
 & 1,  $x$,  $y$, 
$-28 +13 r^2$, $x^2-y^2$, $xy$,
\\
 &
$-(3-r^2) x$, $-(3-r^2) y$,
\\
D2Q13 &
$\frac{1}{12} (202 - 189 r^2 + 35 r^4) x$, 
$\frac{1}{12} (202 - 189 r^2 + 35 r^4) y$, 
\\
 &
$- \frac{1}{2} (280 - 361 r^2 + 154 r^4)$,
$- \frac{1}{12} (65-17r^2)(x^2-y^2)$,
\\
 & 
$- \frac{1}{24}(288 - 1162 r^2 + 819 r^4 - 137 r^6)$
\\
\hline  \end{tabular}
\caption{The orthogonal polynomials for the moments in D2Q$q$ lattice
  Boltzmann models, with $q=1$, 5, 9 and 13. For 3-D cases, see Appendix-2}
\label{tab:D2Q13-polynomials}
\end{table}

Similar expressions can be obtained for 3-D cases (see Appendix~2).
The successive moments can be interpreted as density $\rho$,
components of momentum $\{j_x,j_y\}$, kinetic energy ($E$), components
of the stress tensor, components of heat flux, and so on.

Depending on which situation is to be simulated, we shall consider
that in situations of dimensionality $d$, there are 1, $d+1$ or $d+2$
moments conserved in collisions.  Either $\{\rho\}$, or
$\{\rho,j_x,j_y\}$ or $\{\rho,j_x,j_y,E \}$ allow to simulate respectively
advection--diffusion, athermal Navier--Stokes, Navier--Stokes problems
for $d=2$.  The other moments (non-conserved moments) evolve with
simple linear relaxation:
\begin{equation}
m_p^* = m_p + s_p(m_p^{eq} - m_p)
\label{relmom}
\end{equation}
where $s_p$ is a relaxation rate and $m_p^{eq}$ the equilibrium value of the moment
$m_p$. We consider that $m_p^{eq}$ is a function of the local conserved quantities
and that the relaxation rates are given values.

Numerous papers \cite{DDH92, JKL05, LL00} and practices of numeric
have analyzed the behavior of the model described above in situations
where conserved quantities vary slowly in space and time (on time or
spatial scales large compared to the elementary units $\delta t$ or
$\delta r$.) A popular approach is to follow the kinetic theory
approach with the Chapman--Enskog expansion.
An alternative way proposed by one of us 
performs a Taylor expansion assuming smooth behavior of the conserved
quantities.  

The method involves an expansion of the non-conserved moments in
powers of the time increment $\delta t$ (considered as a small quantity) 
$$m=m^{(0)}+\delta t\ m^{(1)}+\delta t^2\ m^{(2)}+...$$
and to get iteratively the terms $m^{(l)}$. One gets expressions that
involve space derivatives of the conserved moments of increasing order
together with time derivatives. At each step of the process higher
order time derivatives are eliminated by using the results of the previous
step.
%
%

This leads to equivalent PDE's
relating the conserved quantities that are first order in time
derivatives and of desired order in space derivatives (somewhat like
in the hierarchy Euler, Navier--Stokes, Burnett, super--Burnett, {\it
  etc}). A careful analysis of the iterative process allows to state
whether adding more elementary velocities improves the accuracy of the
results already available. Note however that these approaches
(Chapman-Enskog, Taylor expansion, {\it etc}.)  don't give all
the necessary information concerning numerical stability of the method. Useful
results, although not complete, are provided by the study of the
dispersion equation for plane waves summarized in appendix 1.
The equivalent equations method allows to obtain expressions for
higher order terms and thus to discuss resulting inaccuracies and in
some cases ways to improve the models. Some results are presented
below.


\bigskip \bigskip    \noindent {\bf \large    3) \quad  Generation of equivalent equations}  

\smallskip \noindent 
Here we describe the principle of the generation of equivalent equations.

\monitem  {\bf Ingredients} 

To completely define the LBE process, we need the following ingredients~:

-List of elementary velocities, here $\{c_{ix},c_{iy}\}$.
(to simplify writing we take units such that $c_i$ is of the order of 1,
and thus will have just one small parameter to deal with when making expansions.
This is sometimes called the ``acoustic scaling''.)

-Matrix of moments $M$ (of dimension $n\times n$), and $M^{-1}$ its inverse.

-List of moments conserved in collision $W$ (of dimension $n_c$ equal either to 1 or to 3)

-List of equilibrium values of the non-conserved moments $M^{eq}$ ($n-n_c$),
which depend on the local values of the conserved quantities $W$.

-List of relaxation rates for the non-conserved moments $S$ ($n-n_c$).

-Time evolution of the LBE process written in $f$ space as $n$ equations~:
\begin{equation}
  f_j(t+\Delta t, r)=f_j^*(t,r-c_j \Delta t)
\label{LBE}
\end{equation}
where the superscript $^*$ indicates a ``post-collision'' quantity and $\Delta t$
is the small parameter for expansions.

The collision step is performed in moment-space, whereas the propagation step
is performed in f-space.

\monitem  {\bf Iterative process} 

\smallskip \noindent 
We assume ``smoothly varying'' behaviour for all quantities to be dealt with.
%
Then we can expand the relation (\ref{LBE}) at various orders of accuracy 
relative to the small parameter $ \, \Delta t$. 
At order zero, we find that the pre-collision distribution $ \, f\,$ is close to the post-collision
particle distribution  $ \, f^* $:
\begin{equation*}
f^* = f + {\rm O}(\Delta t) \, . 
\end{equation*}
When we re-write this relation in terms of the 
moments $ \, m $, we deduce from the previous relation and the basic iteration of the lattice Bolzmann scheme
\begin{equation}
  m_k^* = m_k + s_k \, ( m_k^{\rm eq} - m_k )  
\label{relaxation}
\end{equation}
the fact that both $ \, m \, $  and $ \, m^* \, $ are close to the 
equilibrium
\begin{equation}
m = m^{\rm eq}  + {\rm O}(\Delta t) \,, \quad  m^* = m^{\rm eq}  + {\rm O}(\Delta t) \,. 
\label{close-to-equilibrium}
\end{equation}
%

\monitem \noindent {\bf Order one} 

\smallskip \noindent 
After this first step, we expand the  relation (\ref{LBE}) at the order one, 
transform the particles into moments  and replace 
the moments in the first order terms by their equilibrium values. We obtain by this way: 
\begin{equation}
 m_k +  \Delta t \, {{\partial  m_k^{\rm eq}}\over{\partial t}}  +  {\rm O}(\Delta t^2) =
 m_k^* - \Delta t \, \sum_{j \ell \alpha}  M_{kj} \, c_j^\alpha \, M^{-1}_{j \ell} \, 
 {{\partial  m^{\rm eq}_\ell}\over{\partial x^\alpha}}   +  {\rm O}(\Delta t^2) \, . 
\label{order-one}
\end{equation}
For the moments that are equilibrium, {\it id est} $ \, m_k \equiv m_k^* $, the relation 
(\ref {order-one}) gives immediatly the equivalent partial differential equations at order one:
\begin{equation}
 {{\partial  m_k^{\rm eq}}\over{\partial t}}  
+  \sum_{\ell \alpha}\big( \sum_j M_{kj} \, c_i^\alpha \, M^{-1}_{j \ell} \big) \, 
 {{\partial  m^{\rm eq}_\ell}\over{\partial x^\alpha}}  =  {\rm O}(\Delta t) \, . 
\label{edp-order-one}
\end{equation}
Moreover, for the moments $ \, m_k \, $ that are not at equilibrium, 
we extract the difference $ \, m_k - m_k^* \,$  from the relations (\ref{relaxation}) and 
(\ref{order-one}). Then a first order expansion for these non-conserved moments emerge:
\begin{equation*} 
m_k =  m_k^{\rm eq} -{{\Delta t}\over{s_k}} \, \Big(  {{\partial  m_k^{\rm eq}}\over{\partial t}} 
+  \sum_{\ell \alpha} \Big( \sum_j M_{kj} \, c_j^\alpha \, M^{-1}_{j \ell} \Big) \, 
 {{\partial  m^{\rm eq}_\ell}\over{\partial x^\alpha}} \Big) +   {\rm O}(\Delta t^2) \, . 
\end{equation*}
It is then usefull to explicit the nonconserved moments afer relaxation, using (\ref{relaxation}) 
and  the previous relation: 
\begin{equation}  \label{m-star}
m_\ell^* =  m_\ell^{\rm eq}   +  \! \Big( 1 - {{1}\over{s_\ell}} \Big)  \,  \Delta t \, 
\Big(  {{\partial  m_\ell^{\rm eq}}\over{\partial t}} 
+  \! \sum_{\ell \beta} \Big( \sum_j M_{\ell j} \, c_j^\beta \, M^{-1}_{j p} \Big) \, 
 {{\partial  m^{\rm eq}_p}\over{\partial x^\beta}} \Big) \! +   {\rm O}(\Delta t^2) \, . 
\end{equation}
%

\monitem \noindent {\bf Expansion at order two and more} 

\smallskip \noindent 
The next step is to expand the relation (\ref{LBE}) up to second order accuracy;
due to (\ref{close-to-equilibrium}), we can replace the moments $ \, m \, $ and $ \, m^* \,$ 
by their equilibrium values for the second order terms. We obtain in this way
\begin{equation} \label{order-two} 
\left \{ \begin{array}{l} \displaystyle 
\!  m_k +  \Delta t \, {{\partial  m_k} \over{\partial t}} 
+  {{\Delta t^2}\over{2}}  \, {{\partial^2  m_k^{\rm eq}}\over{\partial t^2}}  =   m_k^*
- \Delta t \, \sum_{j \ell \alpha}  M_{kj} \, c_j^\alpha \, M^{-1}_{j \ell} \,  {{\partial  m^*_\ell}\over{\partial x^\alpha}} 
 \\ \qquad  \qquad \displaystyle 
+  {{\Delta t^2}\over{2}}  \, \sum_{j \ell \alpha \beta}  M_{kj} \, c_j^\alpha \,  c_j^\beta \,M^{-1}_{j \ell} \, 
 {{\partial^2  m^{\rm eq}_\ell}\over{\partial x^\alpha \, \partial x^\beta}}   +
  {\rm O}(\Delta t^3) \, . 
\end {array} \right .  \end{equation}
In the expansion (\ref{order-two}), there are  three terms of order 2:  
$ \, {{\partial^2  m_k^{\rm eq}}\over{\partial t^2}} \,$  in the left hand side,
the term 
$ \, \sum_{j \ell \alpha \beta}  M_{kj} \, c_j^\alpha \,  c_j^\beta \,M^{-1}_{j \ell} \, 
 {{\partial^2  m^{\rm eq}_\ell}\over{\partial x^\alpha \, \partial x^\beta}} \, $ in the right hand side  
and the term induced by the expansion  (\ref{m-star}) inside the 
first order term  
$  \,  {{\partial  m^*_\ell}\over{\partial x^\alpha}}  \, $ in the left hand side. 
After taking a careful attention of all these terms, we obtain the partial equivalent equations
at order~2. 
For the end of the computation at  second order, we refer to our original contribution \cite{Du08}. 
For the extension at fourth  order in a general nonlinear approach, we refer to \cite{Du22}. 
The extention to linearised schemes at fourth order accuracy has been proposed in  \cite{DL09}.
The algorithm has been simplified in \cite{ADGL14}. In this contribution, we
have used this last version, also called ``Berlin algorithm''. 

\bigskip \noindent 
The basic development is  made in terms of moments~:
\begin{equation}
m_i = m_i^0 + m_i^1 \Delta t + m_i^2 \Delta t^2 + ...
\end{equation}
and we go back and forth between f-space and m-space with matrices $M$ or $M^{-1}$
as necessary.

At order 0, $m_i^0$ is the set of the $n_c$ conserved moments + equilibrium values of
the $n-n_c$ other moments.

We expand Eq.~\ref{LBE} in powers of $\Delta t$ and collect the various powers of $\Delta t$,
The ``propagation'' on the right hand side of Eq.~\ref{LBE} increases the order in $\Delta t$
by one unit, so one gets expressions of the type
\begin{equation}
\sum_p \partial_t^p m_i^{q-p} = \sum A \partial^{p} m_i^{q-1-p}
\label{development}
\end{equation}
where $A$ is an operator expressed in powers of $M^{-1} P  M$ where $P$ is linked to the
velocity set.
This allows to get iteratively the values of the non-conserved moments in terms of space and
time derivatives of the conserved quantities $W$.
There are however unwanted time derivatives of order larger than 1. They are
eliminated iteratively using the results previously derived.
The complexity of the expressions increases very fast with the order of the
iterations, so some care is needed to estimate which contributions can be safely
discarded.
The net result is either 1 or 3 partial differential equations of the $W$ quantities
that are first order in time and high order in space and so can be directly
compared to classic PDE's (Euler, Navier-Stokes, etc...).


\bigskip \bigskip    \noindent {\bf \large    4) \quad  Linear Analysis of Lattice Boltzmann Models in 2D} 
\label{app:stability}

A practical approach to the study of stability is described below.
Several important features of the ability of a LBM model to simulate physical
flows can be obtained for 
 specialized situations that provide a lot of useful information.
Consider a domain with $\{N_x,N_y\}$ active nodes and periodic boundary conditions.
One looks for solutions of the form
\begin{equation}
m_p(i,j,n)=A_p a_p^i b_p^j z^n
\end{equation}
So we take an initial condition periodic in space:
\begin{equation}
m_p(i,j,0)=m_{p0}\phi_{px}^i \phi_{py}^j + M_{p0}
\end{equation}
using phase factors $\phi_{px}=\exp(\imath k_xc_{px})$ and $\phi_{py}=\exp(\imath k_yc_{py})$.
$\{k_x,k_y\}$ can be interpreted as components of the wave vector and $M_{p0}$ 
is linked to a uniform field (say uniform density and constant background velocity,
a situation allowing to test Galilean invariance of the models).

One can compute the moments $m_{p1}$ at time $n=1$. Assuming that the initial amplitudes $m_{p0}$
are small, one linearizes the new values with respect to $m_{p0}$. 
If the components of the wave vector are compatible with the periodicity conditions,
-- $k_x N_x$ and $k_y N_y$ are multiple of $2\pi$ -- then the expressions for the new 
values are the same at all points ($m_{p1}$ within a simple phase factor). The problem
thus simplifies to a $q\times q$ problem.
One gets 
\begin{equation} 
m_{p1} = {E} \, m_{p0} = z \, m_{p0}
\end{equation} 
with a matrix $\tensorfont G$ defined by $q$ equations in ``$f_{\rm space}$'':
\begin{equation}
{G_p}=(I + {M}^{-1}{C} {M}) \phi_{px} \phi_{py}
\end{equation}
and the corresponding one in ``$m_{\rm space}$'',
${E}={M} G {M}^{-1}$. $C$ corresponds to the collision step
and can be obtained from Eq.~(\ref{relmom}).
Under such periodic conditions, analysis can be made at a single node, and so
one just needs to consider $q$-dimensional vectors 
$$\Phi=\{f_0,\cdots,f_{q-1}\}$$
as  elements of phase  space,
together with the scalar product defined as
\begin{equation}
 \langle \Phi^1 | \Phi^2 \rangle 
 =
 \sum_{p=0}^{q-1} f_p^1 f_p^2
\end{equation}

Note that when this is applied to the moments $m_p$, sums of products of small integers
are involved and so there may be accidental degeneracies. It may therefore be quite
useful to determine the rank of the parts of the moment matrix $M$ corresponding
to moments of the same orders.

The determination of the eigenvalues and eigenfunctions of ${E}$ can be done
with the dispersion equation formalism. For particular values of the wave vector, this can
be done analytically. In particular for $k_x=k_y=0$ one gets $z_p=1-s_p$ indicating
that $0\le s_p \le 2$ for stability. For small values of the wave vector, one
can solve the dispersion equation by successive approximations for the roots
$z_l$ close to 1 then compute $\gamma_l=\log(z_l)$ that will be compared to
the predictions of the standard PDE's.
When numerical values of all parameters present in ${E}$ are given, one can
use fast linear algebra packages (for instance in LAPACK) for several values of the
components of the wave vector. Any situation leading to an eigenvalue $z_l$
with modulus greater than 1 is numerically unstable and therefore not suitable for
simulations. It is found that this usually occurs for ``large'' values of $k_x$
or $k_y$ (say between 1 and  $\pi$) so developments in $k_i$ near $k=0$ are often
not able to predict the corresponding instability.


\bigskip \bigskip    \noindent {\bf \large    5) \quad  Analytic results from the linear analysis} 

\monitem \noindent 
{\bf D2Q5 model for advection-diffusion equation in 2D} 
\label{sec:D2Q5}

It has been known for a long time that a 5 velocity (D2Q5) model can be used
to simulate advection-diffusion in 2-D. 
\begin{equation}
  \partial_t \rho + \bm{V} \cdot \bm{\nabla} \rho 
  - \kappa \Delta \rho 
  = 
  0 .
\end{equation}
However it is found that the effective diffusivity $\kappa$ varies as
the square of the advective velocity.  This is not satisfactory so one
can use D2Q9 with adequate expressions for the equilibrium of the
non-conserved moments.

We shall use the following Table~\ref{table-moments-d2q9}, where $V_x$
and $V_y$ are the $x$ and $y$ components of the advective velocity
$\bm{V}$, respectively, and $u$ and $a$ parameters for optimization.


\begin{table} [h]     
\centering
\begin{tabular}{|l|c|c|l|} 
\hline 
Moment & Parity & Rate & \multicolumn{1}{c|}{Equilibrium}
\\
\hline
$\rho$  & $+$ &  0   &   $ \rho$   
\\
$j_x$   & $-$ & $ s_1$ &   $ \rho \, V_x$ 
\\
$j_y$   & $-$ & $ s_1$ &   $ \rho  \,  V_y$ 
\\
$E$     & $+$ & $s_3$ &   
$ \rho  \,  \left( \alpha + 3 V^2 \right)$ 
\\
$p_{xx}$    & $+$ & $s_4$ &  $\rho  \,  ( V_x^2-V_y^2)$ 
\\
$p_{xy}$    & $+$ & $s_4$ &  $\rho \, V_x \,  V_y $ 
\\
$q_x$   & $-$ & $ s_6$ &   $ d_1 \,  \rho  \,  V_x$ 
\\
$q_y$   & $-$ & $ s_6$ &   $ d_1 \,  \rho  \,  V_y$ 
\\
$\varpi$& $+$ & $ s_8$ &   $ \rho \,  
\left( \beta + a\,  V^2 \right)$
\\
\hline  
\end{tabular}
\caption{D2Q9  equilibrium moments  for advection-diffusion, including two parameters,
$u$ and $a$, for further optimization. $V^2 := V_x^2 + V_y^2$.}
\label{table-moments-d2q9} 
\end{table}


The choice of relaxation rates and expressions in terms of velocity
was made in accordance to the symmetry of the set of elementary
velocities (the parity is indicated in the second column to be used for the
particular two-relaxation times (TRT) models).

Applying the Taylor expansion method in the linear case with the
so-called ``Berlin algorithm'' \cite {ADGL14} to third order in space
derivatives and neglecting non linear terms in density, one gets one
equivalent equation for the density:
\begin{equation}
  \partial_t \rho 
  + \bm{V} \cdot \bm{\nabla} \rho
  - \kappa  \,  \Delta \rho
  = 
  {O}(\bm{\nabla}^3 \rho) ,
\label{advdif}
\end{equation}
where  the diffusivity $\kappa$ is independent of velocity:
\begin{equation} 
\kappa=\frac{\alpha+4}{6} \, \left( \frac{1}{s_1}-\frac{1}{2} \right) \, .
\end{equation}  
The next order 
\begin{equation} 
  O(\nabla^3)
  =
  \sum_{\alpha \beta \gamma}  
  H_{\alpha \beta \gamma}(\bm{V})  \, 
  \partial_{\alpha \beta \gamma }^3 \rho ,
  \quad
  \alpha,\, \beta,\, \gamma
  \in \{x,\, y\} ,
\end{equation}  
leads to corrections to advection and thus corresponds to the aim of
the present report.

Considering a plane wave 
$\rho(\bm{r},\, t) = \exp(\gamma t) \exp( \imath  \bm{k} \cdot \bm{r})$
and taking only contributions
linear in velocity in Eq.~(\ref{advdif}), the phase velocity is
\begin{equation} 
\bm{V} \cdot \bm{k} \,  [ 1 + A(\bm{k},\, \hat{\bm{V}}) ] .
\end{equation}  
From now on, we refer to $A(\bm{k},\, \hat{\bm{V}})$ as the
``anomalous advection'' and we try and minimize its magnitude.

This correction factor $A(\bm{k},\, \hat{\bm{V}})$ is a complicated
function depending on the orientations (with respect to the computational grid)
of both the velocity $\bm{V}$
and the wave-vector $\bm{k}$. However it becomes independent of
orientations when
\begin{equation} 
  d_1=-1 
  \quad 
  \mbox{or} 
  \quad 
  \sigma_1 \sigma_4 
  =
  \frac{1}{12} ,
\end{equation}  
where we use the H\'enon parameters \cite{He87}
defined by $ \, \displaystyle \sigma_i= {{1}\over{s_i}}-{1\over2}$.
Note that the condition for isotropy of the shear viscosity of the
standard D2Q9 model leads also to the equivalent value for the
parameter $c_1=-1$.

When $d_1=-1$,
the correction to advection becomes:
\begin{equation} 
A_1
=
\frac{1}{24} \, 
\left[
2+\alpha
+ 
4 (\alpha \sigma_3 - 2 \sigma_4) \sigma_1 
+ 
8 (4 + \alpha) \sigma_1^2 
\right] ,
\end{equation}  
%
and when $\sigma_1  \sigma_4 = \frac{1}{12}$, it reduces to
\begin{equation} 
  A_2
  =
  \frac{1}{72} \, 
  \left[
    7-d_1+3 \, \alpha+12 \, \sigma_1 \, \sigma_3 \, (1-\alpha+d_1)
    - 24 \sigma_1^2 \, (4 + \alpha)
    \right]
  .
\end{equation}  
%
%
Both expressions can be put to 0 by suitable choice of the parameters
$\sigma_i$ provided stability of the process is satisfied.

We mention that the next order in the equivalent equation 
 Eq.~(\ref{advdif}) gives rise to a correction to the viscous
term, allowing to define the ``hyper-diffusivity''. This has been
studied for $V_x=V_y=0$ in ref. \cite{DL09}.
%
The results presented here for the D2Q9 model can be extended to 3-dimensional situations.
The simplest model is based on D3Q7 with elementary velocities $\{0,0,0\}$, 
$\{1,0,0\}$, $\{-1,0,0\}$, $\{0,1,0\}$, $\{0,-1,0\}$, $\{0,0,1\}$, $\{0,0,-1\}$. 
However the effective diffusivity is velocity-dependent. Therefore models based on
D3Q15 or D3Q19 have been proposed. The basic properties of these models for advection-diffusion
and the tuning of parameters to get rid of anomalous advection are summarized in Appendix 2.

\bigskip  \monitem    
{\bf Athermal fluid simulated with D2Q9}
\label{sec:D2Q9}

We start with the common D2Q9 model with 3 conservations defined by
the Table \ref{table-equilibre-d2q9}. 
Applying the Taylor expansion analysis 
up to third order in space derivatives 
leads to a hierarchy of equivalent equations for $\rho, j_x, j_y$
analogous to Equ. \ref{advdif}. 
As we consider only the linear behavior of the three conserved quantities
it is convenient to express the results in terms of matrices for the 
successive orders in space derivatives (shown later as $N_0$, $N_1$, $N_2$ and $N_3$).

\begin{table} [h]
\centering
\begin{tabular}{|c|c|c|c|} 
\hline 
Moment & Parity & Rate & Equilibrium
\\
\hline
$\rho$  & $ + $ &  0   &   $ \rho$   
\\
$j_x$   & $ - $ &  0 &   $ j_x$ 
\\
$j_y$   & $ - $ &  0 &   $ j_y$ 
\\
$E$     & $ + $ & $ s_3$ &   $ \displaystyle \rho \, \bigg( \alpha + 3
\, {{j_x^2+j_y^2}\over{\rho}}  \bigg) $ 
\\
$XX$    & $ + $ & $ s_4$ &  $  \displaystyle {{ j_x^2-j_y^2 }\over
  {\rho}}  $ 
\\
$XY$    & $ + $ & $ s_4$ &  $  \displaystyle {{ j_xj_y }\over { \rho}}
$ 
\\
$q_x$   & $ - $ & $ s_6$ &   $ -j_x$ 
\\
$q_y$   & $ - $ & $ s_6$ &   $ -j_y$ 
\\
$\varpi$& $ + $ & $ s_8$ &   $  \displaystyle \rho \, \bigg( \beta-3
  {{ j_x^2+j_y^2}\over { \rho}}  \bigg) $ 
\\
\hline  
\end{tabular}
\caption{Equilibrium values of  the D2Q9 moments for fluid equations.}  
\label{table-equilibre-d2q9}
\end{table}


The first order, which aims to match Euler's equations, is
\begin{equation} 
  \tensorfont{M}_0 
  + 
  \tensorfont{M}_1
  = 
  \left(
  \begin{array}{ccc}
    \displaystyle 
    \partial_t 
    & 
    \partial_x 
    &
    \partial_y 
    \\
    \displaystyle 
    \frac{(\alpha+4)}{6} \partial_x
    - V_x \bm{V} \!\cdot\! \bm{\nabla}
    & 
    \partial_t + \bm{V} \!\cdot\! \bm{\nabla}
    + V_x  \partial_x
    & 
    V_x \,  \partial_y 
    \\
    \displaystyle 
    \frac{(\alpha+4)}{6} \partial_y
    - V_y \bm{V} \!\cdot\! \bm{\nabla}
    &
    V_y \partial_x 
    & 
    \partial_t
    + \bm{V} \!\cdot\! \bm{\nabla}
    + V_y \partial_y
    \\
  \end{array} 
  \right)
\end{equation} 
higher orders $\tensorfont{M}_l$ are cumbersome and not given here.

Starting from initial conditions
\begin{subequations}
\begin{align}
  \rho(x,y,t) 
  &=
  1 +  \rho_0  
  \exp(\omega t) 
  \cos(\bm{k} \cdot \bm{r}) 
  ,
  \\
  \bm{J}(x,y,t) 
  &=
  \bm{V} 
  + \bm{J}_0
  \exp(\omega t)  \cos(\bm{k} \cdot \bm{r}) 
  ,
\end{align}
\end{subequations}
we apply the matrices $\tensorfont{M}_l$ and show results just for the particular case
where the mean velocity is orthogonal to the wave vector. In addition we
apply a rotation of the axis such that the wave vector is along the axis $Ox$: 
\begin{equation} 
N_0=\left(\begin{array}{ccc}
\omega  & 0 & 0 \cr
0  & \omega & 0 \cr
0 & 0 & \omega \cr
\end{array} \right) \,, 
\end{equation} 
\begin{equation} 
N_1=\left(\begin{array}{ccc}
0  & 1 & 0 \cr
\frac{\alpha +4}{6}  & 0 & 0 \cr
0 & V  & 0 \cr
\end{array} \right)\ k \,, 
\end{equation} 
\begin{equation} 
N_2
=
\left
(\begin{array}{ccc}
0 & 0 & 0 \cr
0 & (\alpha \sigma_3-2\, \sigma_4 ) (1-3 V^2) & 0 \cr
 (\alpha+4) \sigma_4 V  & 0 & -2 \sigma_4  \cr
\end{array} 
\right) 
\frac{k^2}{6}
\, . 
\end{equation} 
Note that Navier-Stokes equations can be expressed just with these three matrices
but without the off-diagonal terms $N_1(3,2)$ and $N_2(3,1)$ and the velocity
in $N_2(2,2)$.

At order 3, taking the usual values of the parameters $\alpha=-2$ and $\beta=1$
in order to simplify the expressions, one gets:
%
\begin{equation} 
N_3
=
\left( \begin{array}{ccc}
\displaystyle 0   & -\frac{1}{18} & 0  \cr
\displaystyle h_0 +h_1\, V^2 &  g_1\, f_1(\theta) & h_3\, V-g_3 \, V \, f_2(\theta) \cr
\displaystyle g_2\, V^2\, f_1(\theta) & h_4 \, V-g_2 \, V^2 f_1(\theta) &g_2\, V\, f_1(\theta) \cr
\end{array} \right)\ k^3
\end{equation} 
with 
\begin{equation}  \left\{ \begin{array}{l}
\displaystyle f_1(\theta)=\sin{4 \,  \theta} \,, \,\,
f_2(\theta)=\sin^2{2 \, \theta} \,,  \\
\displaystyle h_0=\frac{1-3 \, (\sigma_3^2+\sigma_4^2)}{27} \,, \,\,\\  
\displaystyle h_1=\frac{(\sigma_3+3 \, \sigma_4-2 \, \sigma_6)(\sigma_3-\sigma_4)}{6} \,, \,\, 
h_3=\frac{(\sigma_4-2 \,\sigma_6)(\sigma_3-\sigma_4)}{3} \,,  \\ 
\displaystyle g_1=\frac{1- 6 \, \sigma_6 \, (\sigma_3+\sigma_4)}{24} \,, \,\, 
g_2=\frac{1-12 \, \sigma_4 \, \sigma_6}{24} \,, \,\, 
g_3=\frac{1+12 \, \sigma_6 \, (\sigma_3-2 \, \sigma_4)}{24} \,, \\
\displaystyle h_4=\frac{1+6\sigma_4(\sigma_3+\sigma_4-2\sigma_6)}{18}\, \\
\end{array} \right. \end{equation} 
%
%
where the relaxation rates appear as
\begin{equation} 
\sigma_3=\frac{1}{s_{E}}-\frac{1}{2}
\end{equation} 
for the energy mode,
\begin{equation} 
\sigma_4=\frac{1}{s_{XX}}-\frac{1}{2}
\end{equation} 
for components of the stress tensor, and
\begin{equation} 
\sigma_6=\frac{1}{s_{qx}}-\frac{1}{2}
\end{equation} 
for the components of the heat flux.
%
%
This third order matrix becomes independent of the angle $\theta$ for 
\begin{equation}
\sigma_3=\sigma_4\ \quad {\rm and} \quad \sigma_4\ \sigma_6=\frac{1}{12}
\label{isotropy}
\end{equation}
leading to
\begin{equation} 
N_3^{\rm isotropic}
=\left(
\begin{array}{ccc}
0 & 3 (\alpha-2)  & 0 \cr
 (\alpha+4)(\alpha-2) (6\sigma_4^2 - 1)  & 0 & 0 \cr
0 & 6 ( \alpha (12 \sigma_4^2 - 1) -2 )  & 0 \cr
\end{array} 
\right) 
\frac{k^3}{216} 
\end{equation} 
One can then obtain the complex relaxation rate of the waves.

\bigskip 
{\bf Transverse wave and $V$ perpendicular to $k$}

\noindent
At order 1 in $k$, the phase velocity is 0.

\noindent
At order 2 in $k$, the attenuation is $-\sigma_4/3 k^2$, we recover the 
usual shear dynamic viscosity

$$
\nu_0 = \frac{1}{3} \sigma_4  
$$

\noindent
At order 3 in $k$, one gets a phase velocity 
\begin{equation} 
  v_\varphi 
  = 
  \frac{1}{24}
  V \, ( 1-12 \sigma_4 \sigma_6 ) \sin{4 \theta}
  .
\end{equation} 
%

\bigskip 
{\bf Transverse wave and $V$ parallel to $k$}

At order 1 in $k$, the phase velocity is $V$.

At order 2 in $k$, the attenuation corresponds to an effective shear
viscosity
\begin{equation} 
  \nueff
  =
  \frac{1}{3} \sigma_4 \, 
  \left( 1 - 3 V^2 \right) 
  =
  \nu_0 \, \left( 1 - 3 V^2 \right)
  .
\end{equation} 
At order 3 in $k$, the phase velocity is modified (at first order in $V$) by
\begin{equation} 
  \frac{1}{24}
  \left[
    16 \sigma_4 ( \sigma_4-\sigma_6)
    + (1 - 12 \sigma_4 \sigma_6 )
    \, f_2(\theta)  
    \right] 
  V 
  .
\end{equation} 
%
Similarly expressions are readily obtained for acoustic waves 
when $ \, V = 0 \,$ \cite{DL09}; 
cancellation of the corresponding expression occurs for the particular case  
$\nu_{eff} = 1 / \sqrt{108}$, which may be referred to as a ``quartic condition''
which can be seen as the
ancellation of the ``hyper-viscosity''.

\bigskip  \monitem 
{\bf   Athermal fluid simulated with D2Q13}  
\label{sec:D2Q13}

Similar expressions have been derived for the D2Q13
model~\footnote{See Appendix 2 for details on the relaxation step} and
we just give the results of the analysis of the waves, using for the
relaxation rates
\begin{equation} 
\sigma_4=\frac{1}{s_{XX}}-\frac{1}{2}
\end{equation} 
for components of the stress tensor, and
\begin{equation} 
\sigma_6=\frac{1}{s_{qx}}-\frac{1}{2}
\end{equation} 
for the components of the  heat flux, and
\begin{equation} 
\sigma_8=\frac{1}{s_{rx}}-\frac{1}{2}
\end{equation} 
for the components of the ``next''  heat flux.


\bigskip 
{ \bf Transverse wave and $V$ perpendicular to $k$}

At order 1 in $k$, the phase velocity is 0.

At order 2 in $k$, the effective shear viscosity is
\begin{equation} 
  \nueff 
  =
  \frac{3+c_1}{4} \sigma_4 
  \left[
  1 - 
  \frac{12 (7+6q)}{77(3+c_1)} V^2 
  \right]
  =
  \nu_0 
  \left[
  1 - \frac{12 (7+6q)}{77(3+c_1)} V^2 
  \right]
\end{equation} 
showing that one can eliminate the velocity dependence of the effective
shear viscosity for the particular value of the parameter $q=-7/6$.

At order 3 in $k$, there is an additional phase velocity
\begin{equation} 
  v_\varphi 
  =  
  \frac{\sigma_4(89772\, \sigma_6 + 30888\, \sigma_8)-10055}{157080} \, V\, \sin{4\theta}
\end{equation} 
%

\bigskip 
{\bf Transverse wave and $V$ parallel to $k$}

At order 1 in $k$, the phase velocity is $V$.

At order 2 in $k$, the effective shear viscosity is
\begin{equation} 
\nu_{eff}\,=\, \frac{3+c_1}{4} \, \sigma_4 \, \left(1-\frac{12}{77}  \, \frac{7+6q}{3+c_1} \, V^2  \right)
\end{equation} 
showing that the velocity dependence is the same as in the previous case.

At order 3 in $k$, the phase velocity is modified by
\begin{align}
  v_\varphi 
  =   
  &
  \left[ 
    \sigma_4 \, \left( \frac{128-306c_1}{85} \, \sigma_6
    + \frac{306c_1+182}{85}\, \sigma_8 \right)-\frac{31}{102} 
    \right] 
  \,  V \,  f_2(\theta) 
  \nonumber
  \\ 
  &
  + 
  \left[ 
    \sigma_4 \, 
    \left( 
    - \frac{1+c_1}{5} \, \sigma_6
    + \frac{9c_1-1}{20} \,  \sigma_8 
    \right) 
    - \frac{5+3c_1}{48} 
    + \frac{3+c_1}{2}  \,\sigma_4^2 
    \right] 
  \,V \, . 
\end{align}
%



It is possible to remove the angular dependence by taking
\begin{equation} 
\sigma_6=\sigma_8=\frac{1}{12 \, \sigma_4}
\end{equation} 
which leads to an additional phase velocity
\begin{equation} 
v_\varphi \,= \, \frac{3+c_1}{24}  \,  \, (12 \, \sigma_4^2-1)  \,V .
\end{equation} 
The special value $ \, \sigma_4 = {{1}\over{\sqrt{12}}} \, $ allows to get rid of the
additional phase velocity.


In the general case of arbitrary orientations of the wave vector $k$ and of the
advection speed $V$, expressions are quite complicated. Some information
on the relative importance of the corrections to the advection are shown in
Fig.~1. 
The advection term is computed numerically as
\begin{equation} 
  g(\bm{k})  
  = 
  \bm{k} \cdot \bm{V} 
  \left( 1 + h  \, k^2 \right) 
\end{equation} 
with $h$ depending on the orientation of both $\bm{k}$ and
$\bm{V}$. It is represented in Fig.~1 
as solid curve for D2Q13
and a dashed curve for D2Q9 for $\bm{V}$ parallel to $Ox$ and $k$ at angle
$\theta$.

\begin{figure} [H] 
\centerline {\includegraphics[width=0.85 \textwidth, angle=0]{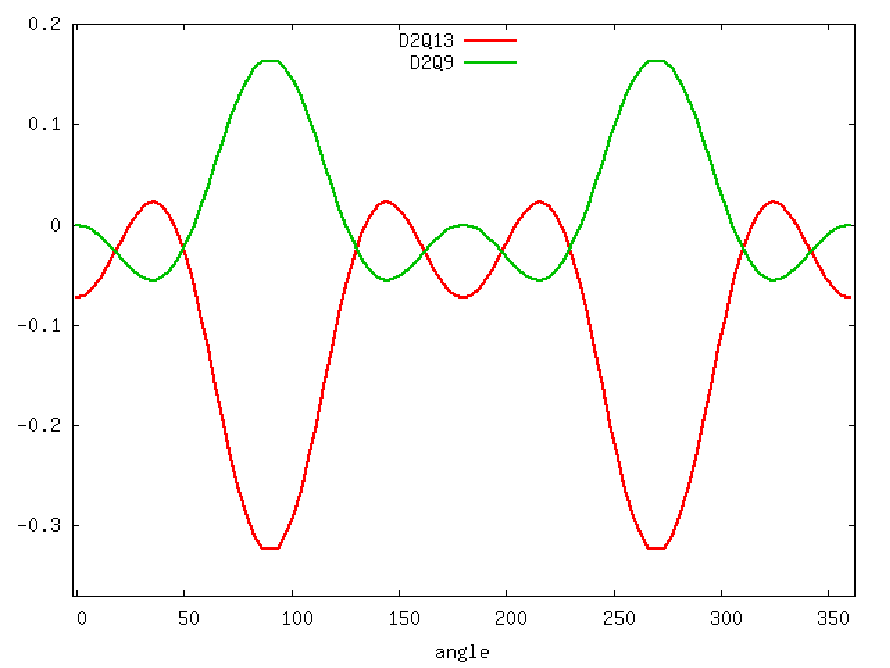}}
\caption{Advection factor for main velocity along $ \, Ox \, $  axis and wave vector
vs angle $\theta$. Dotted line in the absence of anomalous advection. Solid
line contribution $h$ for D2Q13, dashed line for D2Q9.}
\label{gk}  \end{figure}

Some consequences of the correction to advection are presented below.

\endgroup

\begingroup


\bigskip \bigskip   \noindent {\bf \large    6) \quad  Distortion of a Gaussian initial conditions }   

\label{sec:gaussian}

Consider the following Gaussian initial condition
\begin{equation} 
  \Gamma(r,\, 0) 
  = 
  g_0  
  \exp 
  \left[- \left(\frac{r}{r_0} \right)^2 \right] 
  ,
\end{equation} 
centered at the origin $= (0,\, 0)$, where
$r := \sqrt{x^2 + y^2}$ is the distance to the center.
When the Gaussian initial condition $\Gamma(r,\, 0)$ is used as the
initial density $\rho(r,\, 0)$ of the advection-diffusion equation
or the initial stream function $\psi(r,\, 0)$ of the
Navier-Stokes equation, 
 the solution for both cases is
\begin{equation} 
  \Gamma(r,\, t) 
  = 
  g_0
  \frac{r_0^2}{r_0^2 + 4 \chi t}  
  \exp 
  \left[
    - \frac{ (\bm{r} - \bm{V} t) \cdot (\bm{r} - \bm{V} t)}{r_0^2 + 4 \chi t} 
    \right]
\end{equation} 
in the presence of a uniform velocity $\bm{V} := (V_x,\, V_y)$
\cite{LL59}, where $\chi = \kappa$ for the advection-diffusion
equation and $\chi = \nu$ for the Navier-Stokes equation.

The solution $\Gamma(r,\, t)$ is invariant under
rotation.

The results of simulation are shown below for several cases.

\begin{equation} 
  \psi(r,\, 0) 
  = 
  g_0  
  \exp \left[ - \left(\frac{r}{r_0} \right)^2 \right]  
\end{equation} 
or the initial stream function $\psi(r,\, 0)$ for
the Navier-Stokes (centered at the origin \{0,0\}
and $r$ is the distance to the center), evolve as
\begin{equation} 
  \Gamma(r,\, t) 
  = 
  g_0
  \frac{r_0^2}{r_0^2 + 4 \kappa t}  
  \exp 
  \left[
    - \frac{ (\bm{r} - \bm{V} t) \cdot (\bm{r} - \bm{V} t)}{r_0^2 + 4 \kappa t} 
    \right]
\end{equation} 
 or
\begin{equation} 
\psi(r,t) \,= \, g_0 \, \frac{r_0^2}{r_0^2+4  \, \nu  \, t} \, 
\exp \left[ -\frac{(x-V_xt)^2+(y-V_yt)^2}{r_0^2+4 \, \nu  \, t} \right]
\end{equation} 
in the presence of a uniform velocity $\bm{V} := (V_x,\, V_y)$ \cite{LL59}. 
The computed field is invariant by rotation. 
The results of simulation are shown below for several cases.

\newpage 
\bigskip
{\bf Diffuse D2Q9} 

\noindent 
Fig.~\ref{diffud2q9} shows the distribution of $\rho(x,y)$
for three different conditions. The computation is done on a square domain $101^2$
with periodic boundary conditions. The main parameters are:
$\kappa=0.008$, $V_x=0.10,\ V_y=0$ and 3200 time steps. Initial radius is $r_0=5.0$
and initial locations are chosen so that final states do not overlap. The top 
feature is obtained with $q=-1$, the lower feature is obtained with 
$ \, 12 \, \sigma_1 \, \sigma_4=1 \, $
and one can verify that the results are close to rotational invariance.
The right feature
satisfies neither of the isotropy conditions and it is clear that it is not
rotationally invariant.

\begin{figure} [H]  
\centerline {\includegraphics[width=0.85 \textwidth, angle=0]{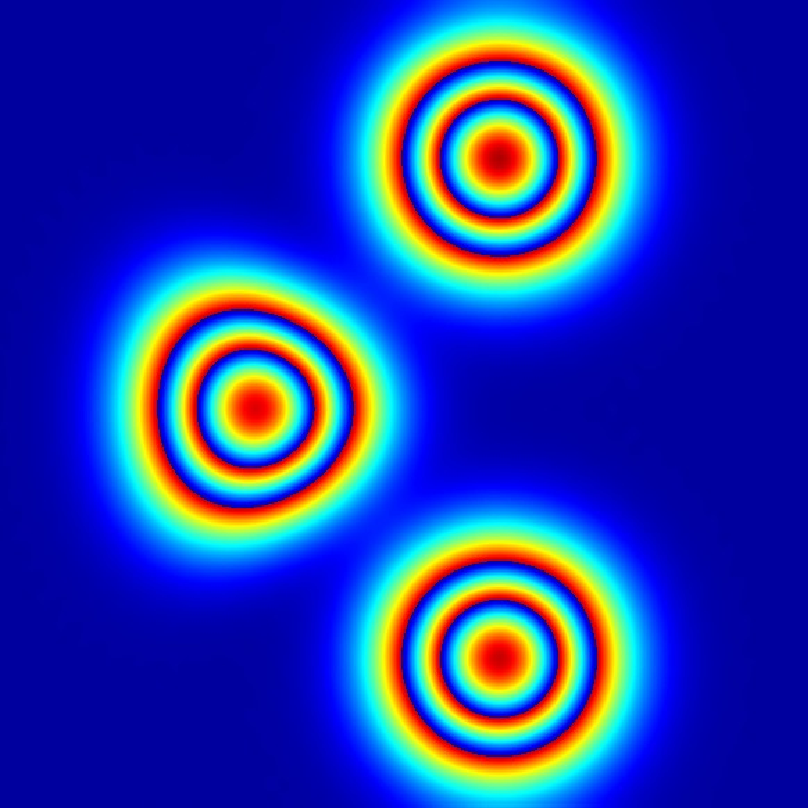}}
\caption{Advection of an initial Gaussian disturbance simulated with diffusive D2Q9
under conditions described in the text. Top and lower features are isotropic
(respectively for $q=-1$ or $12\sigma_1\sigma_4=1$). The middle feature uses
conditions that are not tuned for isotropy.}
\label{diffud2q9}  \end{figure}

\newpage 
\bigskip {\bf Navier--Stokes D2Q9} 

Simulation of the D2Q9 model have been performed in a $301^2$ domain
with periodic boundary conditions. The initial condition is uniform
speed (indicated in the caption), the shear viscosity is $\nu=0.0035$,
the vortex has initial radius $r_0=8.0$. After a number of iterations
the vorticity of the flow is shown in Fig.~\ref{vortd2q9}. The
rotational symmetry is obviously absent when the
condition~\ref{isotropy} is not satisfied, (right feature).  The
feature on the left uses only the second condition of
Eq.\ref{isotropy} as the first one is incompatible with numerical
stability for small shear viscosity.

\begin{figure} [H]  
\centerline {\includegraphics[width=1.0\linewidth]{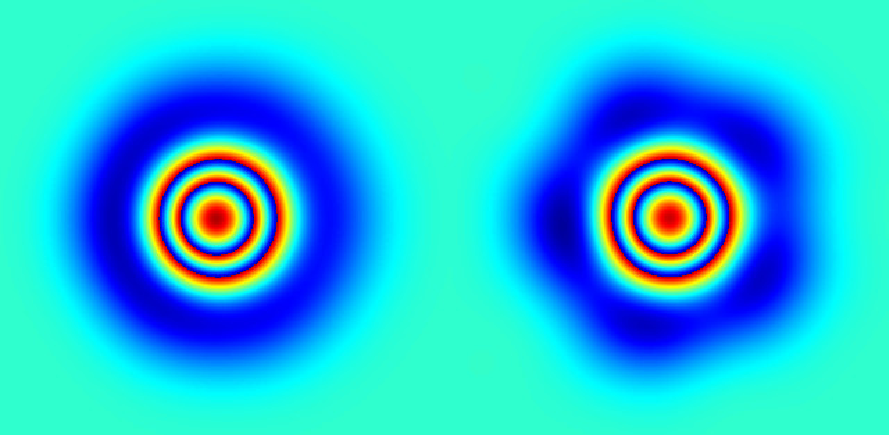}}
\caption{Simulation with D2Q9. Vorticity of the velocity field from an
  initial gaussian stream function after 9000 time steps for an
  advection velocity $\{0.03,0.00\}$. Left with isotropy condition $
  \, 12 \, \sigma_4 \, \sigma_6=1$. Right: arbitrary conditions.}
\label{vortd2q9}  \end{figure}

\bigskip 
{\bf Navier--Stokes D2Q13} 

In a first study, one considers the advection of shear plane waves by a uniform
velocity $V$ parallel to the wave vector. The domain is periodic of size
$240 \times 240$ which corresponds to a smallest wave vector $k_0=2\pi/240$.
Various cases are indicated below with numerical values of the
relative advection either ``experimental''
as determined from simulations or theoretical using expressions given above.

\begin{tabular}{ccccccc}
Case & $k_x/k_0$ &$k_y/k_0$ &$k/k_0$ &  Simulation & Theory & Relative Error \cr
 A & 5 & 12 & 13 & 0.9959 & 0.9960 & 0.01 \% \cr
 B &10 & 24 & 26 & 0.9827 & 0.9840 & 0.13 \% \cr
 C &13 &  0 & 13 & 0.9915 & 0.9917 & 0.02 \% \cr
 D &26 &  0 & 26 & 0.9652 & 0.9666 & 0.15 \% \cr
\end{tabular}

Cases A and B, respectively C and D, correspond to the same orientation of the
wave vector. The data clearly show an increase of the anomaly of the
advection when the wave vector increases and an effect of the orientation.

In a second study, simulation of the D2Q13 model have seen performed 
in a $363^2$ domain with periodic
boundary conditions. The initial condition is uniform speed (indicated in the caption),
the shear viscosity is $\nu=0.003$, the vortex has radius $r_0=11.0$. After a number of
iterations the vorticity of the flow is shown in Fig.~\ref{vortd2q13}. The rotational
symmetry is obviously absent. For comparison the figure also shows what is obtained
without velocity.

\begin{figure} [H] 
\centerline {\includegraphics[width=1.0\linewidth]{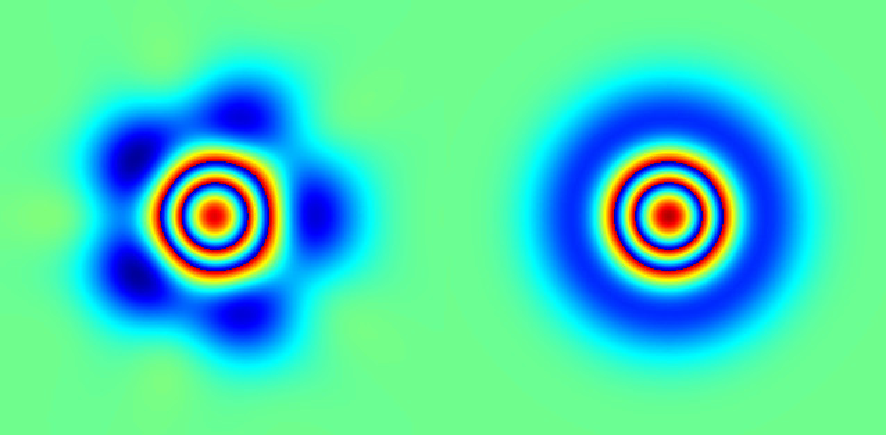}}
\caption{Simulation with D2Q13. Vorticity of the velocity field from
  an initial Gaussian stream function after 2770 time steps. Left with
  an advection velocity $\{0.10,0.00\}$.  Right with no advection.}
\label{vortd2q13} \end{figure}

\bigskip 
{\bf Qualitative interpretation} 
\smallskip \noindent 

To confirm qualitatively the influence of anomalous advection for the present case, 
the advection  is treated in Fourier space. 
The initial stream function $\psi$ can be represented as
\begin{equation} 
\pi \ r_0^2 \ \sum_{k_x,k_y} \, \exp \left[ -r_0^2 \, {{k_x^2+k_y^2}\over{4}}  \right] 
\end{equation} 
and each Fourier component evolves as
\begin{equation} 
\exp \big[ \big(-\nu \, (k_x^2+k_y^2)+ \imath \,  g(k) \, V \big) \, t \big] 
\end{equation} 
For $g(k)$ depending on $k$, the resulting stream function and the
associated vorticity can be computed numerically. An example of such
computations is shown in Fig.~\ref{jouet}.

\begin{figure} [H]
\centerline {\includegraphics[width=1.0\textwidth]{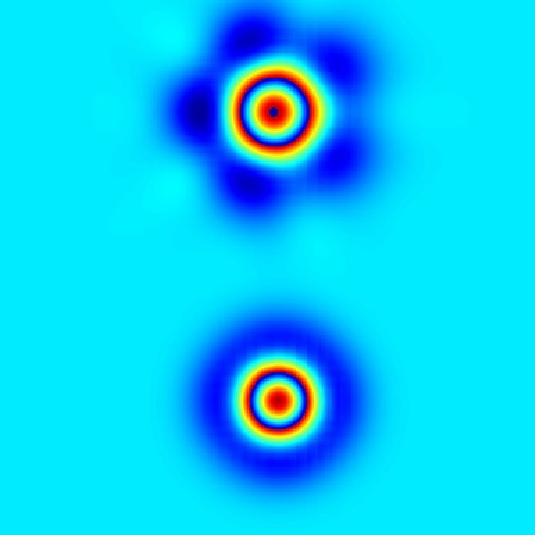}}
\caption{Vorticity of the vortex with main velocity at $14^\circ$ from
  Ox and $r_0=4$ in a domain of size $80\times 80$. Initial state at
  bottom, final state at top. The advection used is $ \, g(k)=1+0.01
  \, \big( \cos (4 \, \theta) - \cos (2 \, \theta) \big)) \, k^2$.}
\label{jouet} \end{figure}

The advection-diffusion case has also been studied in 3-D cases. As recalled earlier,
the simple and popular D3Q7 is inadequate as the diffusivity depends on the square
of the advective velocity, so we give results for D3Q15 and D3Q19 in Appendix 3.

\endgroup

\newpage 
 \bigskip \bigskip   \noindent {\bf \large    Conclusion}   


It has been shown in the present report that lattice Boltzmann models
can be tuned to reduce or in some cases eliminate defects that occur
when they are used to simulate situations of flows with significant
velocities or with features of rather small scales. However the
analysis has been performed only in linearized situations, so that
much work remains to be done for actual nonlinear flows in particular
to estimate the errors due to inaccuracies in the advection which were
pointed by Frisch for the early lattice gas models~\cite{frish}.

\begingroup
\let\clearpage\relax





\bigskip  \bigskip    \noindent {\bf \large    Appendix 1) \quad Moments for the D2Q13 lattice Boltzmann scheme}

\smallskip \noindent 
For the D2Q13 model, we use the moments built with the polynomials
given in Table~\ref{tab:D2Q13-polynomials}.  The equilibrium values
are given in the following Table~\ref{table-moments-d2q13}.

\bigskip   
\begin{table}  [h!]    \centering
\begin{tabular}{|c|c|c|c|} 
\hline 
{\rm Moment} & {\rm Parity} &{\rm Rate} & {\rm Equilibrium} \cr  \hline
$ \displaystyle \rho  $ & $+$  & 0 & $ \displaystyle \rho $ \cr
$ \displaystyle j_x $  & $-$ & 0 &  $ \displaystyle j_x $ \cr
$ \displaystyle j_y $  & $-$  & 0 &  $ \displaystyle j_y $ \cr
E     & $+$  & $ \displaystyle  s_3  $  & $ \displaystyle \alpha\rho+13 \, \frac{j_x^2+j_y^2}{\rho} $ \cr
XX    & $+$  & $ \displaystyle  s_4  $  &  $ \displaystyle \frac{j_x^2-j_y^2}{\rho} $ \cr
XY    & $+$  & $ \displaystyle  s_4  $  &  $ \displaystyle \frac{j_xj_y}{\rho} $ \cr
$ \displaystyle Q_x $   & $-$  & $ \displaystyle s_6 $  & $ \displaystyle 
j_x \, \left(c_1-\frac{36 \, q -35}{77} \, (j_x^2+j_y^2) \right) $ \cr
$ \displaystyle Q_y $   & $-$  & $ \displaystyle s_6 $  & $ \displaystyle 
j_y\, \left( c_1-\frac{36 \, q -35}{77}\,  (j_x^2+j_y^2)\right)  $ \cr
$ \displaystyle R_x $   & $-$  & $ \displaystyle s_8 $   & 
$ \displaystyle j_x\, \left( -\frac{63c_1+65}{24} + q\, j_x^2+\frac{42
  \,q-105}{22}j_y^2 
\right) $ \cr
$ \displaystyle R_y $   & $-$  & $ \displaystyle s_8 $ &  
$ \displaystyle j_y\, \left(  -\frac{63c_1+65}{24}+\frac{42q-105}{22} j_x^2 +q \, j_y^2  \right) $ \qquad \cr
$ \displaystyle E_2 $   & $+$  & $ \displaystyle s_{10} $  & $ \displaystyle \beta \, \rho $ \cr
$ \displaystyle E_3 $   & $+$  & $ \displaystyle s_{11} $  & $ \displaystyle \gamma \, \rho $ \cr
$ \displaystyle XYZ $   & $-$  & $ \displaystyle s_{12} $  & 0 \cr
\hline  \end{tabular}
\caption{Moments of the D2Q13 lattice Boltzmann scheme for fluid flow
including a tuning parameter $q$.}  \label{table-moments-d2q13}
\end{table}

\bigskip \noindent 
The relaxation phase uses the following relaxation rates and equilibrium value, 
such that the speed of sound waves in $ \, c_s=\sqrt {(\alpha+28) / 26 }$
and the shear and bulk viscosities are
\begin{equation} 
  \nu_0 
  = 
  \frac{1}{4} ( c_1+3 )  \sigma_4 ,
  \quad
  \zeta_0
  = 
  \frac{1}{26} (13 c_1 - \alpha+11 )  \sigma_3 
  . 
\end{equation} 
In the presence of a mean velocity, the shear viscosity is
\begin{equation} 
  \nu(V) 
  =
  \nu_0 
  \left[
  1 - 
  \frac{12 (7+6q)}{77 (3+c_1)}  \, V^2 
  \right]
\end{equation} 
leading to optimize the model with $q=-7/6$.


\bigskip  \bigskip    \noindent {\bf \large    Appendix 2) \quad Advection-diffusion  for  three-dimensional situations}
\label{app:a-d-eqn-3d}


\bigskip  {\bf D3Q15 } 

The model follows the usual D3Q15 based of elementary velocities $\{0,0,0\}$, and 
permutations of $\{1,1,1\}$ and of $\{1,0,0\}$.
The moments are computed with the matrix
\begin{equation} 
{M}=
\left(
{\scriptsize
\begin{array}{*{15}{r}}
1&1&1&1&1&1&1&1&1&1&1&1&1&1&1\cr
0&1&$-1$&0&0&0&0&1&$-1$&1&$-1$&1&$-1$&1&$-1$\cr
0&0&0&1&$-1$&0&0&1&1&$-1$&$-1$&1&1&$-1$&$-1$\cr
0&0&0&0&0&1&$-1$&1&1&1&1&$-1$&$-1$&$-1$&$-1$\cr
$-2$&$-1$&$-1$&$-1$&$-1$&$-1$&$-1$&1&1&1&1&1&1&1&1\cr
0&2&2&$-1$&$-1$&$-1$&$-1$&0&0&0&0&0&0&0&0\cr
0&0&0&1&1&$-1$&$-1$&0&0&0&0&0&0&0&0\cr
0&0&0&0&0&0&0&1&$-1$&$-1$&1&1&$-1$&$-1$&1\cr
0&0&0&0&0&0&0&1&1&$-1$&$-1$&$-1$&$-1$&1&1\cr
0&0&0&0&0&0&0&1&$-1$&1&$-1$&$-1$&1&$-1$&1\cr
0&$-4$&4&0&0&0&0&1&$-1$&1&$-1$&1&$-1$&1&$-1$\cr
0&0&0&$-4$&4&0&0&1&1&$-1$&$-1$&1&1&$-1$&$-1$\cr
0&0&0&0&0&$-4$&4&1&1&1&1&$-1$&$-1$&$-1$&$-1$\cr
16&$-4$&$-4$&$-4$&$-4$&$-4$&$-4$&1&1&1&1&1&1&1&1\cr
0&0&0&0&0&0&0&1&$-1$&$-1$&1&$-1$&1&1&$-1$\cr
\end{array} 
}
\right) . 
\end{equation} 
associated to the orthogonal polynomials~:

$$
\begin{array}{|c|c|}
\hline
{\rm Parity } & \cr
\hline
+&1\cr
-&x\cr
-&y\cr
-&z\cr
+&-2+x^2+y^2+z^2\cr
+&2\ x^2-y^2-z^2\cr
+&y^2-z^2\cr
+&x\ y\cr
+&y\  z\cr
+&z\ x\cr
-&x\ (-13/2+5/2\ (x^2+y^2+z^2))\cr
-&y\ (-13/2+5/2\ (x^2+y^2+z^2))\cr
-&z\ (-13/2+5/2\ (x^2+y^2+z^2))\cr
+&16-55/2\ (x^2+y^2+z^2)+15/2\ (x^2+y^2+z^2)^2\cr
-&x\ y\ z\cr
\hline
\end{array}
$$

In the presence of a uniform advective velocity $\{V_x,V_y,V_z\}$, the relaxation
rates $s_i$ and the equilibrium
values of the non-conserved moments are given by the following Table.

\begin{table}[htbp!]    
\centering
\begin{tabular}{|c|c|c|c|} 
\hline
Moment & Parity & Rate & Equilibrium 
\\
\hline
$ \rho $  & $ + $ & $  0  $ & $ 0 $ 
\\
$ j_x   $ & $ - $ & $  s_1  $ & $ \rho \, V_x $ 
\\
$ j_y   $ & $ - $ & $  s_1  $ & $ \rho  \, V_y $ 
\\
$ j_z   $ & $ - $ & $  s_1  $ & $ \rho  \, V_z $ 
\\
$ ee    $ & $ + $ & $  s_5  $ & $ \alpha  \, \rho+\rho  \,
(V_x^2+V_y^2+V_z^2) $ 
\\
$ xx    $ & $ + $ & $  s_6  $ & $ \rho  \, (2  \, V_x^2-V_y^2-V_z^2) $ 
\\
$ yy    $ & $ + $ & $  s_6  $ & $ \rho  \, (V_y^2-V_z^2) $ 
\\
$ xy    $ & $ + $ & $  s_6  $ & $ \rho  \, V_x  \, V_y $ 
\\
$ yz    $ & $ + $ & $  s_6  $ & $ \rho  \, V_y  \, V_z $ 
\\
$ zx    $ & $ + $ & $  s_6  $ & $ \rho  \, V_z  \, V_x $ 
\\
$ q_x   $ & $ - $ & $  s_{11}  $ & $ d_1  \, \rho  \, V_x $ 
\\
$ q_y   $ & $ - $ & $  s_{11}  $ & $ d_1  \, \rho  \, V_y $ 
\\
$ q_z   $ & $ - $ & $  s_{11}  $ & $ d_1  \, \rho  \, V_z $ 
\\
$ d3    $ & $ + $ & $  s_{14}  $ & $ \beta  \, \rho $ 
\\
$ tt    $ & $ - $ & $  s_{15}  $ & $ 0 $ 
\\
\hline  
\end{tabular}
\caption{Equilibrium moments for advective D3Q15.}  
\label{table-6}
\end{table}

This leads to an effective diffusivity
\begin{equation} 
\kappa \, = \, \frac{2+\alpha}{3} \, \sigma_1
\end{equation} 
independent of the velocity.
The analysis of the anomalous advection shows that it can be suppressed for two 
conditions. 

\bigskip  {\bf First case } 
%
%
\begin{equation} 
\sigma_5 \, = \, \frac{4}{(1+3 \alpha)} \, \sigma_6
-\frac{6 \, (2+\alpha)}{(1+3  \, \alpha)} \, \sigma_1
+\frac{3 \, (1+\alpha)}{4 \, (1+3  \, \alpha)}\, {{1}\over{\sigma_1}}
\qquad {\rm for} \quad  d_1=-{7\over3} \, . 
\end{equation} 

\bigskip {\bf  Second case } 
%
%
\begin{equation} 
\sigma_5 \, = \, \frac{10 \, (2+\alpha)}{(3+2\ d_1-5\alpha)} \, \sigma_1
-\frac{15  \, \alpha-2  \, d_1+17}{12 \, (3+2\ d_1-5  \, \alpha)}\, {{1}\over{\sigma_1}} 
\qquad {\rm for} \quad 
\sigma_6 = {{1}\over{12 \, \sigma_1 }} \, . 
\end{equation} 
%

\bigskip  {\bf D3Q19 } 

The model follows the usual D3Q19 based of elementary velocities $\{0,0,0\}$, and 
permutations of $\{1,1,0\}$ and of $\{1,0,0\}$.
The moments are computed with the following matrix ${M}$:
\begin{equation} 
\left( 
{\scriptsize
\begin{array}{*{19}{r}}
  1 &  1 & 1 & 1 & 1 & 1&1&1&1&1&1&1&1&1&1&1&1&1&1
\\
  0&1 & $-1$ & 0 & 0 & 0 & 0 & 1 & $-1$ & 1 & $-1$ & 0 & 0 & 0 & 0 & 1 & 1 &
  $-1$ & $-1$
\\
  0&0 & 0 & 1 & $-1$ & 0 & 0 & 1 & 1 & $-1$ & $-1$ & 1 & $-1$ & 1 & $-1$ & 0 & 0
  & 0 & 0
\\
  0&0 & 0 & 0 & 0 & 1 & $-1$ & 0 & 0 & 0 & 0 & 1 & 1 & $-1$ & $-1$ & 1 & $-1$
  & 1 & $-1$
\\
  $-30$  &$-11$  &$-11$  &$-11$  &$-11$  &$-11$  &$-11$   & 8 & 8 & 8 & 8 & 8 & 8 &
  8 & 8 & 8 & 8 & 8 & 8
\\
  0&2 & 2 & $-1$ & $-1$ & $-1$ & $-1$ & 1 & 1 & 1 & 1 & $-2$ & $-2$ & $-2$ & $-2$ & 1
  & 1 & 1 & 1
\\
  0&0 & 0 & 1 & 1 & $-1$ & $-1$ & 1 & 1 & 1 & 1 & 0 & 0 & 0 & 0 & $-1$ & $-1$
  & $-1$ & $-1$
\\
  0&0 & 0 & 0 & 0 & 0 & 0 & 1 & $-1$ & $-1$ & 1 & 0 & 0 & 0 & 0 & 0 & 0 &
  0 & 0
\\
  0&0 & 0 & 0 & 0 & 0 & 0 & 0 & 0 & 0 & 0 & 1 & $-1$ & $-1$ & 1 & 0 & 0 &
  0 & 0
\\
  0&0 & 0 & 0 & 0 & 0 & 0 & 0 & 0 & 0 & 0 & 0 & 0 & 0 & 0 & 1 & $-1$ &
  $-1$ & 1
\\
  0&$-4$ & 4 & 0 & 0 & 0 & 0 & 1 & $-1$ & 1 & $-1$ & 0 & 0 & 0 & 0 & 1 & 1 &
  $-1$ & $-1$
\\
  0&0 & 0 & $-4$ & 4 & 0 & 0 & 1 & 1 & $-1$ & $-1$ & 1 & $-1$ & 1 & $-1$ & 0 & 0
  & 0 & 0
\\
  0&0 & 0 & 0 & 0 & $-4$ & 4 & 0 & 0 & 0 & 0 & 1 & 1 & $-1$ & $-1$ & 1 & $-1$
  & 1 & $-1$
\\
  0&$-4$ & $-4$ & 2 & 2 & 2 & 2 & 1 & 1 & 1 & 1 & $-2$ & $-2$ & $-2$ & $-2$ & 1 &
  1 & 1 & 1
\\
  0&0 & 0 & $-2$ & $-2$ & 2 & 2 & 1 & 1 & 1 & 1 & 0 & 0 & 0 & 0 & $-1$ & $-1$
  & $-1$ & $-1$
\\
  12&$-4$ & $-4$ & $-4$ & $-4$ & $-4$ & $-4$ & 1 & 1 & 1 & 1 & 1 & 1 & 1 & 1 & 1 &
  1 & 1 & 1
\\
  0&0 & 0 & 0 & 0 & 0 & 0 & 1 & $-1$ & 1 & $-1$ & 0 & 0 & 0 & 0 & $-1$ & $-1$
  & 1 & 1
\\
  0&0 & 0 & 0 & 0 & 0 & 0 & $-1$ & $-1$ & 1 & 1 & 1 & $-1$ & 1 & $-1$ & 0 & 0
  & 0 & 0
\\
  0&0 & 0 & 0 & 0 & 0 & 0 & 0 & 0 & 0 & 0 & $-1$ & $-1$ & 1 & 1 & 1 & $-1$ &
  1 & $-1$
\\
\end{array} 
}
\right). 
\end{equation} 

associated to the orthogonal polynomials~:

$$
\begin{array}{|c|c|}
\hline
{\rm Parity} & \cr
\hline
+&1\cr
-&x\cr
-&y\cr
-&z\cr
+&-30+19\ (x^2+y^2+z^2)\cr
+&2\ x^2-y^2-z^2\cr
+&y^2- z^2\cr
+&x\ y\cr
+&y\ z\cr
+&z\ x\cr
-&x\ (-9+5\ (x^2+y^2+z^2))\cr
-&y\ (-9+5\ (x^2+y^2+z^2))\cr
-&z\ (-9+5\ (x^2+y^2+z^2))\cr
+&(2\ x^2-y^2-z^2)\ (-5+3\ (x^2+y^2+z^2))\cr
+&(y^2-z^2)\ (-5+3\ (x^2+y^2+z^2))\cr
+&12-53/2\ (x^2+y^2+z^2+21/2\ (x^2+y^2+z^2)^2\cr
-&x\ (y^2-z^2)\cr
-&y\ (z^2-x^2)\cr
-&z\ (x^2-y^2)\cr
\hline
\end{array}
$$

In the presence of a uniform advective velocity $\{V_x,V_y,V_z\}$, the relaxation
rates $s_i$ and the equilibrium
values of the non-conserved moments are given by the Table \ref{table-equilibre-d3q19}.

\begin{table}[h!]   
\centering
\begin{tabular}{|c|c|c|l|}  \hline 
Moment  &  Parity &   Rate  &  Equilibrium  \cr
\hline
$ \rho $ & $ + $ & 0 &  $ \rho $ \cr
 $ j_x   $ & $ - $ &  $ s_1  $ &  $ V_x\ \rho $ \cr
 $ j_y   $ & $ - $ &  $ s_1  $ &  $ V_y\ \rho $ \cr
 $ j_z   $ & $ - $ &  $ s_1  $ &  $ V_z\ \rho $ \cr
 $ ee    $ & $ + $ &  $ s_5 $  &  $ \alpha\ \rho+19\ (V_x^2+V_y^2+V_z^2)\ \rho $ \cr
 $ xx    $ & $ + $ & $  s_6  $ & $  (2\ V_x^2-V_y^ 2-V_z^2)\ \rho $ \cr
 $ yy    $ & $ + $ & $  s_6  $ & $  (V_y^2-V_z^2)\ \rho $ \cr
 $ xy    $ & $ + $ & $  s_6  $ & $  V_x\ V_y\ \rho $ \cr
 $ yz    $ & $ + $ & $  s_6  $ & $  V_y\  V_z\ \rho $ \cr
 $ zx    $ & $ + $ & $  s_6  $ & $  V_z\ V_x\ \rho $ \cr
 $ q_x   $ & $ - $ & $  s_{11}  $ & $  d_1\ V_x\ \rho $ \cr
 $ q_y   $ & $ - $ & $  s_{11}  $ & $  d_1\ V_y\ \rho $ \cr
 $ q_z   $ & $ - $ & $  s_{11}  $ & $  d_1\ V_z\ \rho $ \cr
 $ xxe   $ & $ + $ & $  s_{14}  $ & $  0 $ \cr
 $ yye   $ & $ + $ & $  s_{14}  $ & $  0 $ \cr
 $ d3    $ & $ + $ & $  s_{16}  $ & $  \beta\ \rho $ \cr
 $ t_x   $ & $ - $ & $  s_{17}  $ & $  d_2\ V_x\ \rho $ \cr
 $ t_y   $ & $ - $ & $  s_{17}  $ & $  d_2\ V_y\ \rho $ \cr
 $ t_z   $ & $ - $ & $  s_{17}  $ & $  d_2\ V_z\ \rho $ \cr
\hline  \end{tabular}
\caption{Equilibrium moments for the diffusive D3Q19 lattice Boltzmann scheme}  \label{table-equilibre-d3q19}
\end{table}

\noindent 
Applying the same analysis as for D2Q9, one can show that the effective diffusivity
is
\begin{equation} 
\kappa \, = \, \frac{\alpha+30}{57} \,  \left(\frac{1}{s_1}-\frac{1}{2} \right)
\end{equation} 
independent of the velocity $V$.
The order 3 for the equivalent equation includes terms linear in applied velocity
that can be interpreted as corrections to the advection factor.
This correction can be suppressed with two possible sets of parameters.

\bigskip {\bf First case } 

For $d_1=-2/3$ and $d_2=0$, the relaxation rate $s_5$ should satisfy:
\begin{equation} 
\sigma_5 \, = \, \frac{76}{3  \, \alpha-5} \, \sigma_6+\frac{6 \, (\alpha+30)}{5-3 \, \alpha}\sigma_1
+\frac{3 \, (11+\alpha)}{4 \, (3 \, \alpha-5)}\, {{1}\over{\sigma_1}}  
\end{equation} 
where $ \displaystyle \sigma_i = {{1}\over{s_i}}- {1\over2}$ is the H\'enon parameter.

\bigskip {\bf  Second case } 

For $ \, \displaystyle \sigma_6= {{1}\over{12 \, \sigma_1}}$, the relaxation rate $s_5$ should satisfy:
\begin{equation} 
\sigma_5 \,  =  \, \frac{10 \, (\alpha+30)}{21+19 \, d_1-5 \, \alpha} \, \sigma_1-
\frac{279-19 \, d_1+15  \, \alpha}{12 \, (21+19  \, d_1-5  \, \alpha)} \, {{1}\over{\sigma_1}} 
\end{equation} 
Values of the parameters will be constrained by stability conditions, in particular $\sigma_5>0$.

\newpage 
\bigskip {\bf Two Relaxation Times (TRT)   } 

Note that most of the relaxation rates do not appear in the previous conditions, so one can
use the simpler TRT situation (with only two relaxation rates,
one for $+$ parity and one for $-$ parity). The various results shown in this Appendix are summarized
in the table \ref{table-d3q19q15} that applies to the TRT case.

\begin{table}[h!]   
\centering
\begin{tabular}{|l|l|}  
\hline 
Case  & Conditions  
\\ 
\hline 
D3Q19-1   & 
$ \displaystyle d_1 = -\frac{2}{3} \,,\quad   d_2=0   \,,\quad 
\sigma_6 \, =\, \frac{2 \, (30+\alpha)}{27-\alpha}\, \sigma_1
-\frac{11+\alpha}{27-\alpha} \, {{1}\over{4 \, \sigma_1 }}   $ \cr  
D3Q19-2   &  
$ \displaystyle \sigma_1 = {{1}\over{\sqrt{12}}}  \,,\quad  
\sigma_6={{1} \over{\sqrt{12}}}  $\cr  
\cr
D3Q15-1   & 
$ \displaystyle   d_1=-\frac{7}{3}  \,,\quad 
 \sigma_6=2 \, \frac{2+\alpha}{1-\alpha}\sigma_1-\frac{1+\alpha}{1-\alpha}\, {{1}\over{4 \, \sigma_1 }} $ \cr
D3Q15-2   & 
$ \displaystyle  \sigma_1= {{1}\over{\sqrt{12}}}   \,,\quad 
 \sigma_6= {{1}\over{\sqrt{12}}}   $\cr  
\hline  \end{tabular}
\caption{Isotropy of anomalous advection~:~results for the TRT situation.} 
\label{table-d3q19q15}
\end{table}

To be complete, we add some results for the ``hyper-diffusivity'' derived from the equivalent
equations at order 4.

\bigskip  \bigskip    \noindent {\bf \large    Appendix 3)

\hfill   Hyper-diffusivity of the three-dimensional diffusion models}
\label{app:a-d-eqn-3d}

In the absence of an advection velocity, one can easily obtain the 
``hyper-diffusivity'' carrying out the equivalent process to fourth order.
The formula are quite complicated so we only give conditions for
obtaining a null hyper-diffusivity like was done for the shear hyper-viscosity.

\bigskip {\bf D3Q15 } 


%
\begin{equation} 
  \sigma_{11} 
  =
  \frac{(8 \alpha + \beta) + 14 (\alpha + 2) (1 - 6 \sigma_2  \sigma_6 ) }
       {( 8 \alpha + \beta ) (12 \sigma_1 \sigma_6 - 1 )}
       \sigma_1 ,
\end{equation} 
\begin{align} 
  \sigma_5 
  =
  &
  \frac{\sigma_1}{4 \left[ 5 \, (\alpha + 2) (3 \alpha+1) 
    (1 - 12 \sigma_1 \sigma_6 ) + 30 (\alpha + 1)  + 2 ( \beta - 1)
      \right]}
  \times
  \nonumber
  \\
  &
  \left\{ 
  60 ( \alpha + 2)^2 (12 \sigma_1 \sigma_6 - 1) \sigma_1^2  
  - 960 ( \alpha + 2)  \sigma_1^2 \sigma_6^2  
  \right.
  \nonumber
  \\
  &
  \left.
  + 4 ( 2 \beta - 40 \alpha + 68 - 45 \alpha^2) 
  \sigma_1 \sigma_6
  + 15  (\alpha + 2)  \alpha 
  \right\}
\end{align}
%


\bigskip  {\bf D3Q19 } 

\begin{equation} 
  \sigma_{11}
  =
  - \frac{1}{19} 
  \sigma_1
  {{ 84 (\alpha + 30)  \sigma_1 \sigma_6 - 95 \beta - 52 \, \alpha-420 }
\over{ ( 2 \alpha + 5 \beta )  ( 12 \sigma_1 \sigma_6 - 1 ) }} 
\end{equation} 
\begin{align}
  \sigma_5
  = 
  &
  -\frac{1}{4 
    \left\{ 
      84 \, (\alpha + 30) (3 \alpha-5)  \sigma_1 \sigma_6 
      - 21 
      \, \alpha^2 - 722 \, \beta-937 \, \alpha-546 
      \right\}
    \sigma_1 } 
  \times
  \nonumber
  \\
  &
  (1008 \, (\alpha + 30 )^2  \sigma_1^3 \, \sigma_6 
  - 84 (  (\alpha + 30)^2 + 304 ( \alpha + 30)  \sigma_6^2 )
  \sigma_1^2 
  \nonumber
  \\ 
  &
  + 4 (32676  - 63 \alpha^2 - 512  \alpha + 722 \beta) \, \sigma_6
  \, \sigma_1 
  + 21 (\alpha + 30) (\alpha -8)) 
\end{align}


These expressions can be simplified for the TRT case. One obtains the same results
for the two models:
\begin{equation} 
\sigma_1=\sigma_{11}=\frac{1}{\sqrt{12}} \,,\quad 
\sigma_5=\sigma_6=\frac{1}{\sqrt{3}} \,,\quad 
\end{equation} 
Note that one gets the same value of $\sigma_1$ as in Table 1, but a different one for
$\sigma_6$. It is thus not possible to have at the same time no anomalous convection and no
hyper-diffusivity.


\bigskip \bigskip      \noindent {\bf  \large  References }



\endgroup

\end{document}